\documentclass[12pt,twoside]{amsart}
\usepackage{amssymb} 
\usepackage{color}              
\usepackage{epsfig}

\nonstopmode

\numberwithin{equation}{section}
\font\tengothic=eufm10 scaled\magstep 1
\font\sevengothic=eufm7 scaled\magstep 1 
\newfam\gothicfam
\textfont\gothicfam=\tengothic 
\scriptfont\gothicfam=\sevengothic

\newtheorem{theorem}{Theorem}[section]
\newtheorem{lemma}[theorem]{Lemma}
\newtheorem{proposition}[theorem]{Proposition}
\newtheorem{corollary}[theorem]{Corollary}
\newtheorem{conjecture}[theorem]{Conjecture}

\theoremstyle{definition} 
\newtheorem{remark}[theorem]{Remark}
\newtheorem{example}[theorem]{Example}

\newcommand {\RR}{\mathbb{R}}
 
\newcommand {\zz}{\mathbb{Z}} 
\newcommand{\kk}{\mathbb{K}}
\newcommand  {\NN}{\mathbb{N}} 
\newcommand {\PP}{\mathbb{P}}

\newcommand{\cut}{{\mathrm{Cut}^{\square}}}
\newcommand{\ca}{\mathcal{A}}

\begin{document}

\title{Toric Geometry of Cuts and Splits}

\author{Bernd Sturmfels}
\address{Department of Mathematics, University of California, Berkeley, CA 94720}
\email{bernd@math.berkeley.edu}
\author{Seth Sullivant}
\address{Society of Fellows and Department of Mathematics, Harvard University,
Cambridge, MA 02138}
\email{seths@math.harvard.edu}

\begin{abstract}
Associated to any graph is a toric ideal whose generators record relations among the cuts of the graph.  We study these ideals and the geometry of the corresponding toric varieties.  
Our theorems and conjectures  relate the
combinatorial structure of the graph and the corresponding cut polytope 
 to algebraic properties of the ideal.  
 Cut ideals generalize toric ideals arising in phylogenetics and the study of
 contingency tables. 
\end{abstract}
\maketitle

\section{Introduction}

With any finite graph $G = (V, E)$ we associate a projective toric variety
$X_G$ over a field $\kk$ as follows. The coordinates $q_{A|B}$ of the ambient
projective space are indexed by unordered partitions
$A|B$ of the vertex set $V$. The dense torus has
two coordinates $(s_{ij},t_{ij})$ for each edge $\{i,j\} \in E$.
The polynomial rings in these two sets of unknowns are 
\begin{eqnarray*} &
\kk[q]\, \quad := \quad \,\kk\bigl[ 
\,q_{A|B} \, \, | \, \, A \cup B = V, A \cap B =\emptyset \,\bigr], \\ &
\kk[s,t] \quad := \quad \kk\bigl[ \,s_{ij},\, t_{ij} \, \, |
 \, \, \{i,j\} \in E \,\bigr]. \qquad \qquad
\end{eqnarray*}
Each partition $A|B$ of the vertex set $V$ defines a subset 
${\rm Cut}(A|B)$ of the edge set $E$.
Namely, ${\rm Cut}(A|B)$ is the set of edges
$\{i,j\}$ such that $i \in A, j \in B$ or $i \in B,j\in A$.
The variety we wish to study is specified by the following homomorphism of
polynomial rings:
\begin{equation}
\label{monomialmap}
\phi_G \,: \,\kk[q] \,\rightarrow \,\kk[s,t] \, , \qquad
q_{A|B} \,\,\,\, \mapsto
 \prod_{\{i,j\} \in {\rm Cut}(A|B)} \!\!\!\!\! s_{ij}
\,\,\, 
\cdot \!\!\! \prod_{\{i,j\} \in E \backslash {\rm Cut}(A|B)}
 \!\!\!\!\!\! t_{ij} . 
\end{equation}
One may wish to think of $s_{\cdot \cdot}$ and $t_{\cdot \cdot}$ as abbreviations for
``{\bf s}eparated'' and ``{\bf t}ogether''.
The kernel of $\phi_G$ is a homogeneous toric ideal $I_G$
which we call the {\em cut ideal} of the graph $G$.
We are interested in the projective toric variety $X_G$
which is defined  by the cut ideal $I_G$.

\begin{example}
\label{K4}
Let $G = K_4$ be the complete graph on four nodes, so
$V = \{1,2,3,4\}$ and $E = \bigl\{12,13,14,23,24,34 \bigr\}$. The ring map $\phi_{K_4}$ is specified by
$$
\begin{matrix}
                   q_{|1234} \, \mapsto \,    t_{12} t_{13} t_{14} t_{23} t_{24} t_{34},\,
&\quad &  q_{1|234} \, \mapsto \,    s_{12} s_{13} s_{14} t_{23} t_{24} t_{34},
\\                 q_{12|34} \, \mapsto \,    t_{12} s_{13} s_{14} s_{23} s_{24} t_{34},
& \quad &  q_{2|134} \, \mapsto \,    s_{12} t_{13} t_{14} s_{23} s_{24} t_{34},
\\                  q_{13|24} \, \mapsto \,    s_{12} t_{13} s_{14} s_{23} t_{24} s_{34},
& \quad &  q_{3|124} \, \mapsto \,    t_{12} s_{13} t_{14} s_{23} t_{24} s_{34},
\\                q_{14|23} \, \mapsto \,    s_{12} s_{13} t_{14} t_{23} s_{24} s_{34} ,
&\quad &  q_{4|123} \, \mapsto \,    t_{12} t_{13} s_{14} t_{23} s_{24} s_{34} .\end{matrix}
$$
The cut ideal for the complete graph on four nodes is the principal ideal
$$ I_{K_4} \quad = \quad \left<   \,
q_{ | 1234}\, q_{12|34} \, q_{13|24} \, q_{14|23} \,
- \, q_{1|234} \, q_{2| 134} \, q_{3| 124} \, q_{123|4} \,\right>.$$
Thus the toric variety $X_{K_4}$ defined by 
$I_{K_4}$ is  a quartic hypersurface in $\PP^7$. \qed
\end{example}

\begin{example}
\label{C4}
Let $G = C_4$ be the 4-cycle with edges $E= \{12,23,34,14\}$.  
The ring map $\phi_{C_4}$ is derived from $\phi_{K_4}$
in Example \ref{K4} by setting $s_{13} = t_{13} = s_{24} = t_{24} = 1$, and we find
$$I_{C_4} \,\,\,= \,\,\, \left<  \,
q_{ |1234}\, q_{13|24} - q_{1|234} \,q_{124|3}, \,
q_{ |1234} \,q_{13|24} - q_{123|4}  \,q_{134|2} ,\, 
q_{ |1234} \,q_{13|24} - q_{12|34}  \,q_{14|23} \, \right> .$$
Thus the toric variety $X_{C_4}$ is a complete intersection of three quadrics in $\PP^7$. \qed
\end{example}

We usually take the vertex set $V$ of our graph $G$ to be $\,[n] := \{1,2,\ldots,n\}$,
so that $\kk[q]$ is a polynomial ring in $2^{n-1}$ unknowns, and
$\kk[s,t]$ is a polynomial ring in $\,2|E| \leq n(n-1)\,$ unknowns.
Each edge $\{i,j\} \in E$ corresponds to a projective line
$\PP^1$ with homogeneous coordinates $(s_{ij}:t_{ij})$,
and the ring map $\,\phi_G \,$ represents a 
rational map from the product of projective lines
$\,(\PP^1)^{|E|}\,$ into the high-dimensional projective space $\,\PP^{2^{n-1}-1}$.
The image of this map is our toric variety $X_G$, which has dimension 
$\, |E| \leq n(n-1)/2\, $ in $\,\PP^{2^{n-1}-1}$.

The algebraic properties of its ideal $I_G$ and the
 geometry of $X_G$ are determined by
the {\em cut polytope} $\cut(G)$, which is  the convex hull 
in $\RR^{|E|}$  of the {\em cut semimetrics} $\,\delta_{A|B} $. Here
 $A|B$ runs over all unordered partitions of $V$, and
   $\,\delta_{A|B} \in \{0,1\}^{|E|}\,$ is defined as follows:
$$ \delta_{A|B}(\{i,j\}) = 1 \,\,\,\hbox{if}\, \,\,|A \cap \{i,j\}| =1 \quad\, \hbox{and} \,\quad
 \delta_{A|B}(\{i,j\}) = 0 \,\,\,\hbox{ otherwise.} $$
Indeed, the convex hull  of the exponent vectors 
in $\phi_G$ is affinely isomorphic to $\cut(G)$.
In Example \ref{K4} and \ref{C4}, we find that  $\cut(K_4)$ is the cyclic
$6$-polytope with $8$ vertices, and $\cut(C_4)$ is the
$4$-dimensional crosspolytope (which is the dual to the $4$-cube).

The cut polytope $\cut(G)$ is  well-studied  in 
combinatorial optimization, and is a central player in the book
{\em Geometry of Cuts and Metrics} by D\'eza and Laurent \cite{Deza1997}.
The title of this paper is a reference to their book, and reflects our desire to
import this body of work into commutative algebra and algebraic statistics.
In particular, we explore  the extent to which the known polyhedral structure of $\cut(G)$ 
can be used to determine algebraic results about the cut ideals $I_G$.
For instance, the known fact that $\cut(G)$ is full-dimensional implies that $\dim X_G = |E|$.
A more significant example of such an algebraic result is
derived from recent work of the second author \cite{Sullivant2006a}:

\begin{theorem} 
The cut ideal $I_G$ has a squarefree reverse lexicographic initial ideal if and only if 
the graph $G$ is free of $K_5$ minors and every induced cycle in $G$ has length three or 
four.  In this case, every reverse lexicographic initial ideal of $I_G$ is squarefree.
\end{theorem}

\begin{proof}
The initial ideal of a toric ideal is squarefree if and only if the corresponding regular triangulation of the associated polytope is unimodular \cite[\S 8]{Sturmfels1995}.
 Since the symmetry group of $\cut(G)$ is transitive on its vertices, 
 the cut polytope $\cut(G)$ has a unimodular revlex (pulling) triangulation if and only if every revlex triangulation is unimodular \cite[Cor. 2.5]{Sullivant2006a}.  A polytope all of whose revlex triangulations 
are unimodular is called {\em compressed}. Now simply apply the
classification of  compressed cut polytopes given in \cite[Thm. 3.2]{Sullivant2006a}.
\end{proof}

The outline of this paper is as follows.  In Section 2, we describe how 
generating sets (Markov bases) and Gr\"obner bases of the cut ideal $I_G$ 
can be computed when the graph $G$ admits a certain clique-sum decomposition.
The key tool here is the toric fiber product which was introduced in \cite{Sullivant2006b}.  
  In Section \ref{sec:C}, we summarize the results of our computational experiments,
  and we outline some conjectures which were suggested by our computations.

In  the last two sections we present applications to  algebraic statistics.
In Section \ref{sec:graph} we relate cut ideals to the binary graph models of \cite{Develin2003}
and to Markov random fields.
In Section \ref{sec:split} we relate cut ideals to phylogenetic models on split systems,
due to Bryant \cite{Bryant2005}. These generalize the binary Jukes-Cantor
 models which were studied in \cite{BW} and  \cite{Sturmfels2005}.


\section{Clique Sums and Toric Fiber Products}

Our goal in this section is to relate the graph-theoretic operation of taking clique sums to the ideal-theoretic operation of taking the toric fiber product, as explained in \cite{Sullivant2006b}.  This operation will serve as a tool for reducing the computation of the cut ideals $I_G$ 
to cut ideals of smaller graphs (and that, hence, involve fewer  indeterminates). 

Let $G_1 = (V_1,E_1)$ and $G_2 = (V_2,E_2)$ be graphs such that 
$V_1 \cap V_2 $ is a clique of both graphs.
 The new graph $G = G_1 \# G_2$ with vertex set $V = V_1 \cup V_2$ and edge set 
 $E = E_1 \cup E_2$ is called the {\em clique sum} of $G_1$ and $G_2$ along 
 $V_1 \cap V_2$.   If the cardinality of $V_1 \cap V_2$ is 
$k{+}1$, this operation is 
 also called a {\em $k$-sum}  of the graphs.  We suppose throughout that $k \leq 2$.

We now explain how binomials in the cut ideal $I_G$ can be constructed
from binomials in the smaller ideals  $I_{G_1}$ and $I_{G_2}$.  Consider an arbitrary binomial of degree $d$ in the
first smaller cut ideal $I_{G_1}$, say
$$ {\bf f} \quad = \quad \prod_{i = 1}^d q_{A_i |B_i} \, -\, \prod_{i = 1}^d q_{C_i | D_i}. $$
Since $V_1 \cap V_2$ is a clique in $G_1$ of cardinality $\leq 3$, we can permute the
unknowns and partitions so that
 $A_i \cap V_1 \cap V_2   = C_i \cap V_1 \cap V_2$ for all $i$.  This is a consequence of the fact that $I_{K_{k+1}}$ is the zero ideal for $k \leq 2$.
 For any ordered list $EF$ of $d$ partitions of $V_2 \backslash V_1$,
 $$ EF \quad = \quad \bigl(  E_1 | F_1 , \,  E_2 | F_2 , \, \ldots, \, E_d| F_d \bigr), $$
 we define a new binomial which is easily seen to be in 
 the cut ideal $I_G$ of the big graph:
 $$   {\bf f}^{EF} \quad := \quad
  \prod_{i = 1}^d q_{A_i  \cup E_i |B_i \cup F_i}  \,-\, \prod_{i = 1}^d q_{C_i \cup E_i | D_i \cup F_i}.$$
This construction works verbatim if we switch the 
components $G_1$ and $G_2$, so that,
for any binomial ${\bf f}$ in $I_{G_2}$ and any
ordered list $EF$ of $\,{\rm deg} ({\bf f})\,$ partitions of $V_1 \backslash V_2$, we get a binomial
${\bf f}^{EF}$ in $I_G$. Moreover, if 
 ${\bf F}$ is any set of binomials in $I_{G_1}$ or in $I_{G_2}$
 then we define
 \begin{equation}
 \label{liftdef}
 {\rm Lift}({\bf F}) \quad := \quad \bigl\{\,
  {\bf f}^{EF} \, \, | \, \, {\bf f} \in {\bf F}\,, \,\,EF = \{E_i | F_i \}_{i = 1}^{\deg {\bf f} }  \, \bigr\}
 \end{equation}
as the union of all binomials of the form ${\bf f}^{EF}$ described above.

We also define an additional set  ${\rm Quad}(G_1,G_2)$ of
 quadratic binomials  in $I_G$ as follows.
Let  $A|B$ be any unordered partition of $V_1 \cap V_2$, let
$ C_1|D_1$ and $E_1|F_1$ be any ordered partitions of $V_1\backslash V_2$,
and  let $ C_2|D_2$ and $E_2|F_2$ be any ordered partitions of $V_2 \backslash V_1$. Then
\begin{equation}
\label{quadric}
 q_{A \cup C_1 \cup C_2 | B \cup D_1 \cup D_2} \cdot 
      q_{A \cup E_1 \cup E_2 | B \cup F_1 \cup F_2 } \, - \,
     q_{A \cup E_1 \cup C_2 | B \cup F_1 \cup D_2} \cdot 
      q_{A \cup C_1 \cup E_2 | B \cup D_1 \cup F_2 }      
       \end{equation}
 is in ${\rm Quad}(G_1,G_2)$, and these are all the binomials in ${\rm Quad}(G_1,G_2)$.
For each fixed $A|B$, we can express the quadrics
(\ref{quadric}) as the $2 \times 2$-minors of a certain matrix
$\,(q_{\bullet |  \bullet })\,$ of
format $2^{|V_2 \backslash V_1 |} \times 2^{|V_1\backslash V_2|}$.
 The following theorem will be our main result in Section~2.

 \begin{theorem}\label{thm:cutfiber}
Let $G = G_1 \# G_2$ be a $0$, $1$,  or $2$-sum of $G_1$ and $G_2$ and suppose that 
${\bf F}_1$ and ${\bf F}_2$ are binomial generating sets for
the smaller cut ideals $I_{G_1}$ and $I_{G_2}$.
Then 
$$ {\bf M} \quad = \quad {\rm Lift}({\bf F}_1)
\,  \cup\, {\rm Lift}({\bf F}_2 )\, \cup \,{\rm Quad}(G_1, G_2)$$
is a generating set for the big cut ideal $I_G$.  Furthermore, if ${\bf F}_1$ and ${\bf F}_2$ are Gr\"obner bases for $I_{G_1}$ and $I_{G_2}$ then there exists a term order such that 
${\bf M}$ is a Gr\"obner basis for $I_G$.
\end{theorem}

\begin{remark}
If the intersection graph  $G_1 \cap G_2$ is not a clique of cardinality $\leq 3$, then
 it is generally not possible to lift every binomial in $I_{G_1}$ and $I_{G_2}$ to the cut ideal $I_G$.
\end{remark}

Before presenting the proof of Theorem \ref{thm:cutfiber}
 we discuss several examples and corollaries.

\begin{example}\label{ex:segre}
If $G = G_1 \# G_2$ is a zero sum,  then its cut ideal $I_G $
is the usual Segre product of $I_{G_1}$ and $I_{G_2}$. Indeed,
in this case the singleton $V_1 \cap V_2$ has only one ordered partition 
and ${\rm Quad}(G_1,G_2)$ is the ideal of $2 \times 2$-minors
of the corresponding matrix $(q_{\bullet|\bullet})$.
 For instance,
if $G_1$ is the graph with one edge $\{1,2\}$ and
$G_2$ is the graph with one edge $\{2,3\}$, so $V_1 \cap V_2 = \{2\}$,
then $\, I_G = \langle{\rm Quad}(G_1,G_2) \rangle $\, is generated by the determinant of
$$ (q_{\bullet|\bullet}) \quad = \quad \begin{pmatrix}
 q_{|123} & q_{1|23} \\
 q_{12|3} & q_{2|13}  \end{pmatrix}. $$
Now suppose that $G$ is any tree with $n$ leaves.
Iterating the zero sum construction from $n = 3$ to $n > 3$,
we see that $X_G$ is 
  the Segre embedding of $\,(\PP^1)^{n-1}\,$ into 
$\,\PP^{2^{n-1}-1}$. \qed
\end{example}

Further generalization of Example \ref{ex:segre} leads to the following result.

\begin{corollary}
The toric variety $X_G$ is smooth if and only if $G$ is free of $C_4$ minors.
\end{corollary}

\begin{proof}
We first prove the if-direction. If $G$ is free of $C_4$ minors, so all its simple
cycles have length three, then $G$ can be built from $K_2$ and $K_3$
by taking repeated $0$-sums. Both the ideals $I_{K_2}$ and $I_{K_3}$ are zero and live
in polynomial rings with two and four unknowns respectively.
Thus $X_{K_2}$ is $\PP^1$ and $X_{K_3}$ is $\PP^3$.
The $0$-sum construction amounts to taking Segre products, hence
$$ X_G \quad = \quad \PP^1 \times \PP^1 \times \cdots \times \PP^1 \,\times \,
                                         \PP^3 \times \PP^3 \times \cdots \times \PP^3 .$$
This Segre variety is smooth. The only-if direction says that any smooth
$X_G$ has this special form. To prove this, 
suppose that $G$ has $C_4$ as a minor.  Then either $G$ has an induced cycle of length
 $n \geq 4$, or $G$ has as an induced subgraph the complete graph $K_4$ or the graph 
 which is obtained from $K_4$ by removing one edge.
 Let $H$ denote this induced subgraph.  
Using a forward reference to Lemma \ref{lemma:induced},
we note that $\cut(H)$ is a face of $\cut(G)$.
Therefore, it suffices to prove that $X_H$ is not smooth.
We saw in the Introduction that $\cut(K_4)$ and
$\cut(C_4)$ are not simple. Using the familiar
characterization of toric singularities \cite[\S 2.1]{Fulton}, this implies that
the corresponding toric varieties $X_H$ are not smooth.
The same can be checked for cycles of length $n \geq 5$. 

In the remaining case,  $H = K_4 \backslash \{14\}$ is the
$1$-sum of the triangle on $\{1,2,3\}$ and the triangle on $\{2,3,4\}$.
Its variety $X_H$ is the complete intersection of two quadrics in $\PP^7$:
$$ I_H \quad  = \quad \left< \,
{\rm det}
\begin{pmatrix}
q_{|1234}  & q_{1|234} \\
q_{4|123} & q_{14|23} 
\end{pmatrix}
\, , \,\,
{\rm det}
\begin{pmatrix}
q_{2|134}  & q_{12|34} \\
q_{13|24} & q_{3|124} 
\end{pmatrix}\,
\right>.
$$
The singular locus of $X_H$ consists of the two
$3$-planes in $\PP^7$ where these matrices are zero.
The cut polytope $\cut(H)$ is the free join of two squares, a
non-simple $5$-polytope.
\end{proof}

The following example naturally generalizes the graph 
$H = K_4 \backslash \{14\}$ we just discussed.

\begin{example}\label{ex:k5minuse}
Let $G = K_5 \backslash \{15\}$ be the graph on five vertices obtained from 
the complete graph by deleting an edge.  Thus $G$ is the $2$-sum
of the complete graph $G_1$ on $V_1 = \{1,2,3,4\}$ and the complete
graph $G_2$ on $V_2 = \{2,3,4,5\}$. Since $I_{K_4}$ is generated by a quartic, we deduce that $I_G$ is generated by quadrics and quartics. There are four quadrics:
\begin{eqnarray*} {\rm Quad}(G_1,G_2) \,\,= &
\bigl\{
q_{15|234} \, q_{|12345} \, - \, q_{1|2345} \, q_{5|1234}\, , \,\,\,
q_{34|125} \,  q_{2|1345} \, - \, q_{12|345} \, q_{25|134}\\ & \quad \,
q_{24|135} \,   q_{3|1245} \, - \, q_{13|245} \,  q_{35|124} \, , \,\,\,
q_{23|145} \,  q_{4|1235} \, - \, q_{14|235}  \,  q_{45|123} \bigr\}.
\end{eqnarray*}
The ideals $I_{G_1}$ and $I_{G_2}$ are each generated
by a single quartic, as in Example \ref{K4}, and 
${\bf F}_1$ and ${\bf F}_2$ are the singletons consisting
of these quartics. Now, the  set $\, V_2 \backslash V_1 = \{5\}\,$
has two ordered partitions, namely $5|$ and $|5$, so there
are $2^4 = 16$ ordered lists of ordered partitions $E|F$. Each
defines a quartic in $I_G$, so 
${\rm Lift}({\bf F}_1)$ consists of $16$ quartics, such as
\begin{eqnarray*}
& {\bf f}_1 \quad = \quad         
  q_{|12345} \, q_{34|125} \, q_{24|135} \, q_{23|145} \,-\,
  q_{1|2345} \, q_{25|134} \, q_{35|124} \, q_{45|123} ,\\
& {\bf f}_2 \quad = \quad  
   q_{5|1234} \, q_{12|345} \, q_{13|245} \, q_{14|235} \, - \,
    q_{15|234} \, q_{2|1345} \, q_{3|1245} \, q_{4|1235}.
    \end{eqnarray*}
Likewise, ${\rm Lift}({\bf F}_2)$ consists of $16$ quartics, and these include
\begin{eqnarray*}
& {\bf f}_3 \quad = \quad        
         q_{1|2345} \,q_{25|134}\, q_{35|124}\, q_{45|123} \, - \,
         q_{15|234} \, q_{2|1345} \, q_{3|1245} \, q_{4|1235}, \\
& {\bf f}_4 \quad = \quad
         q_{|12345} \, q_{34|125} \, q_{24|135} \,  q_{23|145} \, - \,
         q_{5|1234} \, q_{12|345} \,  q_{13|245} \, q_{14|235}.
 \end{eqnarray*}
We conclude that the set ${\bf M}$ in Theorem \ref{thm:cutfiber} 
consists of $36$ binomials, and these binomials generate $I_G$.
However, they are not a minimal generating set, because of the relation
$$ {\bf f}_1\, -\, {\bf f}_2 \,+\, {\bf f}_3\, -\, {\bf f}_4 \quad = \quad 0 .$$
The set of $35$ binomials obtained by
removing any of the ${\bf f}_i$ is a 
minimal generating set for $I_G$.
We also find that the minimal free resolution of $I_G$ has the following Betti diagram:
\begin{verbatim}
     total: 1 35 134 200 134 35 1
         0: 1  .   .   .   .  . .
         1: .  4   .   .   .  . .
         2: .  .   6   .   .  . .
         3: . 31 128 200 128 31 .                     (Macaulay 2 output).
         4: .  .   .   .   6  . .
         5: .  .   .   .   .  4 .
         6: .  .   .   .   .  . 1
\end{verbatim}
Thus the toric $9$-fold $X_G \subset \PP^{15}$ is
arithmetically Gorenstein. The degree of $X_G$ is $80$.
\qed
 \end{example} 

 \smallskip

 We shall derive Theorem \ref{thm:cutfiber} from the  results in \cite{Sullivant2006b}.  
 Specifically, we shall identify the cut ideal of
 $G = G_1 \# G_2$ as a toric fiber product.
  We begin by reviewing the set-up of \cite{Sullivant2006b}.
 Let $r > 0$ be a positive integer and $s, t \in \NN^r$ be two vectors of positive integers.  Let
$$\kk[x] \,=\, \kk\bigl[\,x^i_j \, \, | \, \, i \in [r], j \in [s_i]\, \bigr]
\quad  \mbox{ and }  \quad \kk[y] \,=\, \kk \bigl[\,y^i_k \, \, | \, \, i \in [r], k \in [t_i]\, \bigr]$$
be two polynomial rings which have a compatible
$d$-dimensional multigrading
$$ \qquad 
\deg(x^i_j) \,=\, \deg(y^i_k) \,=\, {\bf a}^i \,\,\,\in \,\,\,\zz^d
\qquad \qquad \mbox{(for $i=1,2,\ldots,r$).} $$
We abbreviate the collection of degree vectors
by $\,\ca = \{{\bf a}^1,{\bf a}^2, \ldots, {\bf a}^r \} \subset \zz^d$.
	
If $I$ and $J$ are homogeneous ideals of $\kk[x]$ and $\kk[y]$ respectively, then
the quotient rings $R = \kk[x]/I$ and $S = \kk[y]/J$ are also multigraded by $\ca$.
Consider the polynomial ring
$$\kk[z] \,\,\,= \,\,\,\kk \bigl[\,z^i_{jk}  \, \, | \, \, i \in [r], j \in [s_i], k \in [t_i] \,\bigr] ,$$
and consider the $\kk$-algebra homomorphism
$$\phi_{I,J} \,:\, \kk[z] \rightarrow R \otimes_\kk S\,, \,\,\,\,
  z^i_{jk}  \mapsto  x^i_j \otimes y^i_k. $$
  The kernel of $\phi_{I,J}$ is called the {\em  toric fiber product} of $I$ and $J$
  and is denoted
$$I \times_\ca J \quad = \quad \ker(\phi_{I,J}).$$
The following statement combines Theorem 2.8 and Corollary 2.10 in \cite{Sullivant2006b}.

\begin{theorem}  \label{journalofalgebra}
Suppose that the set $\ca$  of degree vectors
is linearly independent.  Let ${\bf F}_1$ be a homogeneous generating set for $I$ and 
${\bf F}_2$ be a homogeneous generating set for $J$.  Then
$$ {\bf M} \,\,\, = \,\,\,{\rm Lift}({\bf F}_1) \,\cup \,{\rm Lift }({\bf F}_2) \,\cup \,{\rm Quad}_\ca $$
is a homogeneous generating set for $I \times_\ca J$.  Furthermore, 
if $\,{\bf F}_1$ and ${\bf F}_2$ are Gr\"obner bases for $I$ and $J$, 
then there exists a term order such that 
${\bf M}$ is a Gr\"obner basis for $I \times_\ca J$.
\end{theorem}

Here ${\rm Quad}_\ca$ is the collection of quadrics 
$\, z^i_{jk} z^i_{lm} - z^i_{jm} z^i_{lk} \,$ which generates
$\,\left<0 \right> \times_\ca \left<0\right> $. The sets ${\rm Lift}({\bf F}_i)$ have a nice
description in terms of tableaux which is given in \cite[\S 2]{Sullivant2006b}.

\medskip

 \noindent {\sl Proof of Theorem \ref{thm:cutfiber}. }
 Suppose $G = G_1 \# G_2$ with vertex set $V = V_1 \cup V_2$ and edge set 
 $E = E_1 \cup E_2$ where $V_1 \cap V_2$ is a clique of size $k+1$ in both graphs.
We set $d = \binom{k+1}{2}+1$ and $r = 2^{k-1}$, and we define
$\ca$ as the vector configuration corresponding to the
vertices of the cut polytope $\cut(K_{k+1})$ of the clique.
The $\ca$-grading on $\kk[q]$ is defined by restricting 
the product in (\ref{monomialmap}) to those
edges $\{i,j\}$ which lie in $E_1 \cap E_2$.
In other words, the degree of $q_{A|B}$ is the
vertex of $\cut(K_{k+1})$ which is indexed by the
partition $\,A \cap V_1 \cap V_2|B \cap V_1 \cap V_2$.

 The configuration $\ca$ of degree vectors is linearly independent if and only
 if the cut polytope $\cut(K_{k+1})$ is a simplex if and only if $k \leq 2$.
Theorem \ref{journalofalgebra} requires the
set $\ca$ to be linear independent. This explains the crucial
  hypothesis $k \leq 2$ in Theorem \ref{thm:cutfiber}.

All three cut ideals $I_G$, $I_{G_1}$ and $I_{G_2}$ are 
homogeneous with respect to the indicated grading.  
We will show that $I_G$ is the toric fiber product of $I_{G_1}$
and $I_{G_2}$, in symbols,
\begin{equation}
\label{isToricFP}
 I_G \,\,\,= \,\,\, I_{G_1} \times_\ca I_{G_2}.
\end{equation}
Let $A_1 | B_1$ and $A_2| B_2$ be partitions of $V_1 $ and $V_2$
such that $\deg(q_{A_1 | B_1} )  = \deg(q_{A_2 | B_2})$. 
 Since $V_1 \cap V_2$ is connected, this implies (possibly after
relabeling) that $A_1 \cap V_1 \cap V_2 = A_2 \cap V_1 \cap V_2$.
This means that $A|B$ with $A = A_1 \cup A_2$ and $B = B_1 \cup B_2$ 
is a partition of $V$, and we have
\begin{equation}
\label{magiceqn}
\phi_{G_1}(q_{A_1|B_1}) \cdot \phi_{G_2}(q_{A_2|B_2}) \quad = \quad
\phi_{G}(q_{A|B}) \cdot \phi_{G_1 \cap G_2}(q_{
A_1 \cap V_1 \cap V_2|B_1 \cap V_1 \cap V_2}).
\end{equation}
This is an identity of monomials in the polynomial ring
$\kk[s,t]$ associated with the big graph $G$, and it
is verified by plugging in the definition of
the monomial map $\phi_\bullet$
in (\ref{monomialmap}).

The ring map which defines 
the toric fiber product $\, I_{G_1} \times_\ca I_{G_2}\,$
can be written as follows:
$$ \phi_{I_{G_1},I_{G_2}}\,: \,
 \kk[q] \, \rightarrow \, \kk[s,t] \, , \,\,\,\,
q_{A|B} \, \mapsto \, \phi_{G_1}(q_{A_1|B_1}) \cdot 
\phi_{G_2}(q_{A_2|B_2}) . $$
Since (\ref{magiceqn}) holds and since
$\phi_{G_1 \cap G_2}(q_{A_1 \cap V_1 \cap V_2|B_1 \cap V_1 \cap V_2})\,$
divides $\,\phi_{G}(q_{A|B})$, the
unknowns $s_{ij}$ or $t_{ij}$
with $\{i,j\} \in E_1 \cap E_2$ can appear in
$\,\phi_{G_1}(q_{A_1|B_1}) \cdot \phi_{G_2}(q_{A_2|B_2}) \,$
only with exponent $2$.
If we replace these unknowns $s_{ij},t_{ij}$
by their square roots in the monomial map
$\, \phi_{I_{G_1},I_{G_2}}\,$ then the kernel
remains unchanged, and we get the monomial map
$\, \phi_{G} : \kk[q] \, \rightarrow \, \kk[s,t] $.
We conclude that
$\, {\rm ker}( \phi_{G}) \, = \,
 {\rm ker}( \phi_{I_{G_1},I_{G_2}})$, which is
our claim (\ref{isToricFP}).
Since  the configuration $\ca$ is linearly independent,
we have thus derived Theorem \ref{thm:cutfiber}
from Theorem \ref{journalofalgebra}.
\qed

\medskip

The proof of Theorem  \ref{journalofalgebra}
given in \cite{Sullivant2006b} reveals the possible
choices of term orders which create a Gr\"obner basis
for $I_G$ from given Gr\"obner bases  ${\bf F}_1$  
of $I_{G_1}$ and  ${\bf F}_2$ of $I_{G_2}$.
First of all, the passage from a binomial ${\bf f}$
in ${\bf F}_i$ to the corresponding binomials
${\bf f}^{FE}$ in ${\rm Lift}({\bf F}_i)$
is compatible with the  choice of leading terms,
that is, we  declare the leading term of 
${\bf f}^{FE}$ to be the one coming from the
leading term of ${\bf f}$. In this manner we
specify a family of partial term orders  on $\kk[q]$.
We then choose any tie-breaking term order on $\kk[q]$
which makes the set  ${\rm Quad}(G_1,G_2)$ into 
a Gr\"obner basis. Since these quadrics are the
$2 \times 2$-minors of matrices $(q_{\bullet | \bullet})$
whose entries are disjoint sets of unknowns there
are many such choices of term orders. Now, the
term order on $\kk[q]$ which is gotten by refining 
the partial term order by the tie-breaker has
the desired property that ${\bf M}$ is a
Gr\"obner basis for $I_G$.


\section{Computations and Conjectures}\label{sec:C}

Upon encountering a new family of ideals, our first instinct is to use computer algebra  to gain a better ``feel'' for the way the structure of the ideals depends on the parameters defining the ideals.  The  parameter for 
the cut ideal $I_G$ is the graph $G$, and we are interested in how the combinatorial structure of $G$ determines the algebraic structure of $I_G$.  To this end, we 
undertook an exploration of the cut ideals by computing generating sets, Gr\"obner bases, free resolutions, and normalizations, using the programs {\tt 4ti2} \cite{Hemmecke2003}, {\tt CoCoA} \cite{CoCoA}, {\tt Macaulay 2} \cite{Grayson}, and {\tt Normaliz} \cite{Bruns}.  In this section, we summarize the results of our computations,
and we offer a number of conjectures that arise from looking at the resulting data.

\subsection{Computations}

The results are summarized in Table 1 below.  The first column lists the graphs which we analyzed.  These were all graphs on six or fewer vertices that are not clique-sum decomposable with a clique of size $\leq 3$.  The notation of the form $G_k$ comes from
the {\em Atlas of Graphs} \cite{Read1998}.  However, if a graph has a more standard shorthand, we preferred to use the more easily identifiable abbreviations.  The notations we used are:
\begin{itemize}
\item $K_{l}$  $\qquad\, $ Complete graph, 
\item $K_{l_1,...,l_m}$ $\,$ Complete $m$-partite graph,
\item $C_l$  $\qquad \,\,$ Cycle of length $l$,
\item $\widehat{G}$  $\qquad \,\,$ Suspension of $G$ over a point,
\item $G \times H$  $\,\,$ Cartesian product graph.
\end{itemize}

\noindent The columns in the table list the following features of the cut ideal $I_G$:
\begin{itemize}
\item[\bf 2-6]  Number of minimal generators of $I_G$ in degrees 2, 4, 6, 8, and 10.
\item[\bf 7]  $ \mu(I_G)\, = \,\,$ Largest degree of a minimal generator of $I_G$.
\item[\bf 8]  Codimension (height) of $I_G$.
\item[\bf 9]  Projective dimension of $I_G$.
\item[\bf 10]  Degree (multiplicity) of $I_G$.
\item[\bf 11]  Whether the semigroup algebra $\kk[q]/I_G$ is normal.
\item[\bf 12]  Whether the semigroup algebra $\kk[q]/I_G$ is Cohen-Macaulay.
\item[\bf 13]  Whether the semigroup algebra $\kk[q]/I_G$ is Gorenstein.
\end{itemize}

Blank spots in the table are entries that we were unable to compute.

\begin{table}
\bigskip
\label{thetable}
\begin{tabular}{|c||c|c|c|c|c|c|c|c|c|c|c|c|} 
\hline
   & 2 & 4 & 6 & 8 & 10 & $\mu$ & codim & pdim & deg & nor & CM & Gor \\
\hline
\hline
$K_3$ & 0 & 0 &  0 & 0 & 0 & 0 & 0 & 0  & 1  & Y  & Y  & Y \\
\hline
\hline
$C_4$ & 3 & 0 & 0 & 0 & 0 & 2 & 3 & 3  & 8  & Y  & Y  & Y  \\
\hline
$K_4$ & 0 & 1 & 0 & 0 & 0 & 4 & 1 & 1  & 4  & Y  & Y  & Y \\
\hline
\hline
$C_5$ & 30 &0  &0 & 0&  0& 2 & 10 & 10   & 52  & Y  & Y  & N   \\
\hline
$K_{2,3}$  & 19  & 0 & 0  & 0  & 0 & 2 &  9 & 9  & 72  & Y  & Y  & Y \\
\hline
$G_{48}$ & 14  & 4   & 0   & 0   & 0  & 4 & 8 & 8  & 60  & Y  & Y  & N \\
\hline 
$\widehat{C_4}$ & 8  & 8  & 0   & 0    & 0  & 4 & 7 & 7  & 64  & Y  & Y  & N \\
\hline 
$K_5$  & 0  & 20  & 40  & 0  & 0  & 6 & 5 & 15  & 128  & N  & N  &  N   \\
\hline 
\hline
$C_6$  &  195 & 0  & 0  & 0  & 0  & 2 & 25  & 25   & 344   & Y  & Y  & N  \\
\hline 
$G_{129}$ & 146  & 0   & 0  & 0   & 0   & 2   & 24 & 24  & 712  & Y  & Y  &  N\\
\hline 
$K_{2,4}$ &  111 & 0   & 0   & 0   & 0   & 2 & 23 & 23  & 1152   & Y  & Y  &  Y \\
\hline 
$G_{151}$ & 118  & 16  & 0  & 0  & 0  & 4 & 23 & 23  & 912  & Y  & Y  & N \\
\hline 
$G_{153}$ & 132  & 12  & 0  & 0  & 0  & 4 & 23 & 23  &  608 & Y  & Y  &  N \\
\hline 
$G_{154}$ &  111 & 16  & 0  & 0  & 0  & 4 & 23  & 23  & 1280   & Y  & Y  & Y  \\
\hline 
$G_{170}$ & 94  &  64  & 0  & 0   & 0  & 4 & 22 &  22 & 1344  & Y  & Y  & N \\
\hline 
$G_{171}$ &  100 & 28  & 0  & 0  & 0  & 4 & 22 & 22  & 976   & Y  & Y  & N \\
\hline 
$G_{173}$ & 90  & 52  & 0  & 0  & 0  & 4 & 22 & 22  & 1440  &  Y &  Y &  N \\
\hline 
$K_2 \times K_3$ &  90 & 52   & 0  & 0  & 0  & 4 & 22 & 22   & 1440  & Y  & Y  & N \\
\hline 
$K_{3,3}$ & 63  & 72  & 0   & 0   & 0   & 4 & 22 & 22  & 3168  & Y  & Y  & Y  \\
\hline 
$G_{186}$ &  72  & 196   & 0  & 0  & 0  & 4 & 21 &  21 & 1984  & Y  & Y  &  N \\
\hline 
$\widehat{C_5}$ &  80 & 40   & 0  & 0  & 0  & 4 & 21 & 21   &  1232 & Y  & Y  &  N \\
\hline 
$G_{188}$ & 64  & 114   &  0  & 0  & 0  & 4 & 21 & 21  & 1856  & Y  & Y  & N \\
\hline 
$G_{189}$ & 54  & 246  & 0  & 0  & 0  & 4 & 21 & 21  & 2976  &  Y & Y  & N \\
\hline 
$G_{190}$ & 76  & 128   & 0  & 0   & 0  & 4 & 21 &  21 & 1600  & Y  & Y  & N \\
\hline 
$G_{194}$ & 60  &  207 & 160  & 0  & 0  & 6 & 20 &   & 3184  & N  & N  & N \\
\hline 
$\widehat{ K_{2,3} }$ & 44  & 420  & 0  & 0  & 0  & 4 & 20 & 20  & 3360  & Y  & Y  & N \\
\hline 
$G_{198}$ &  48 & 336  & 0  & 0  & 0  & 4 & 20 &  20 & 3040  & Y  & Y  & N \\
\hline 
$G_{199}$ & 44  & 337  & 80  & 0  & 0  & 6 & 20 &   & 3760  & N  & N  & N \\
\hline 
$G_{203}$ & 32  & 473  & 160  & 0   & 0  & 6  & 19 &   & 5696   & N  &  N &  N \\
\hline 
$K_{2,2,2}$ &  24 & 1096  & 0  & 0  & 0  & 4 & 19 & 19   & 6144  &Y   & Y   &  N \\
\hline 
$G_{206}$ & 16  & 671  & 320  & 0  & 0  & 6 & 18 &   & 11520  &  N &  N & N \\
\hline 
$G_{207}$ & 8  & 436  & 2872  & 0  & 0  & 6 & 17 &   & 23104  & N  &  N & N \\
\hline 
$K_6$ & 0  &  260 &  3952  & 846  & 480  & 10 & 16 &   & 52448  & N  &  N & N \\
\hline 
\end{tabular}
\bigskip
\caption{Algebraic properties of cut ideals $I_G$ for graphs $G$
with up to six vertices.}
\end{table}
\medskip

If $G$ is a small clique-sum decomposable graph then we can break it into pieces
that are listed in Table 1. This tells us the degrees of the minimal generators of 
cut ideal $I_G$, but it does not tell all invariants of $I_G$.
To be precise, although Theorem \ref{thm:cutfiber} shows that
$$ {\bf M} \quad =  \quad {\rm Lift}({\bf F}_1)
\,  \cup\, {\rm Lift}({\bf F}_2 )\, \cup \,{\rm Quad}(G_1, G_2)$$
generates the cut ideal $I_G$, when ${\bf F}_1$ and ${\bf F}_2$ are minimal generating sets of $I_{G_1}$ and $I_{G_2}$, the set ${\bf M}$ need not generate minimally.
This happens in Example \ref{ex:k5minuse}.  Furthermore, we do not know how taking toric fiber products affects the Cohen-Macaulay type.  For instance, the usual Segre product of two Gorenstein ideals need not be Gorenstein.

\subsection{Conjectures}

We now present some conjectures inspired by our computations.  Our main observation is that many of the coarse invariants of the cut ideals seem to be preserved under taking minors of the underlying graph.  Recall that a graph $H$ is a minor of $G$ if $H$ can be obtained from $G$ by deleting and contracting edges.  By the Robertson-Seymour Theorem on graph minors \cite{Robertson2004}, we may hope to characterize the class of graphs whose cut ideals satisfy some algebraic property by a finite list of excluded minors.  

The protypical example of such a conjecture concerns the maximal degree of a binomial appearing in 
a minimal generating set of the cut ideal $I_G$. This number is $\mu(I_G)$.

\begin{conjecture}\label{conj:minors}
The set of graphs $G$ such that $\mu(I_G) \leq k$ is minor-closed for any $k$.
\end{conjecture}

As evidence for Conjecture \ref{conj:minors}, note that two operations related to taking graph minors 
amount to taking faces of the corresponding cut polytopes.

\begin{lemma} \label{lemma:induced}
\begin{enumerate}
\item   If $H$ is an induced subgraph of $G$ then $\cut(H)$ is a face of $\cut(G)$.
\item  If $H$ is obtained from $G$ by contracting an edge then $\cut(H)$ is a face of $\cut(G)$.
\end{enumerate}
\end{lemma}

\begin{proof}
For part (2), intersect $\cut(G)$ with the hyperplane $x_{ij} = 0$ where $ij$ is the contracted edge.
For part (1), intersect $\cut(G)$ with the hyperplanes $x_{ij} = 0$ for all edges $ij$ in $G$ not incident to $H$, together with one extra condition $x_{ij} = 0$ for each connected component of $G \setminus H$, where $ij$ is an edge incident to said connected component and $H$.  
\end{proof}

This implies that generating degrees can only go down when
 passing to an induced subgraph or when contracting an edge:

\begin{corollary} \label{corcont}
\begin{enumerate}
\item  If $H$ is an induced subgraph of $G$ then $\mu(I_H) \leq \mu(I_G)$.
\item  If $H$ is obtained from $G$ by contracting an edge then  $\mu(I_H) \leq \mu(I_G)$.
\end{enumerate}
\end{corollary}

\begin{proof}
For any two toric ideals, we 
always have the inequality $\mu(I_{\mathcal{B}}) \leq \mu(I_{\mathcal{A}})$ whenever $\mathcal{B}$ is a face of $\ca$.  Thus, the desired inequalities are a direct consequence of Lemma \ref{lemma:induced}.
\end{proof}

Therefore, to prove Conjecture \ref{conj:minors}, it would suffice to show that 
generating degrees are nonincreasing upon the deletion of edges.  Note that the face property does not hold when deleting an edge, as seen by comparing Examples \ref{K4} and \ref{C4}.

\begin{conjecture} \label{conj:delete}
Let $H$ be obtained from $G$ by deleting an edge.  Then $\mu(I_H) \leq \mu(I_G)$.
\end{conjecture}

The smallest instance of Conjecture \ref{conj:minors}, namely
$k=2$ concerns those graphs $G$ whose cut ideal $I_G$ is generated
by quadrics. We propose the following simple characterization:

\begin{conjecture}\label{conj:degree2}
The cut ideal $I_G$ is generated by quadrics
 if and only if $G$ is free of $K_4$ minors (i.e. if and only if $G$ is series-parallel).
\end{conjecture}

If a graph $G$ has $K_n$ as a minor, then that minor can be realized 
by a sequence of edge contractions only.  
By Corollary \ref{corcont} (2), the cut ideal of
every graph with a $K_4$ minor has a minimal generator 
of degree $4$.  Thus, to prove Conjecture \ref{conj:degree2} 
we must show that graphs without $K_4$ minors have quadratically
generated cut ideals.
 Graphs free of $K_4$ minors are known as \emph{series-parallel graphs}.
Every series-parallel graph can be built from $K_2$ by successive series 
and parallel extensions.  The series extensions are just $0$-sums.
Hence, to prove Conjecture \ref{conj:degree2}, it would suffice to show 
that $\mu(I_G)$ does not increase when performing a parallel extension. 

Another conjecture, along the same
lines as Conjecture \ref{conj:degree2},
concerns quartic generators. 

\begin{conjecture} \label{conj:K_5}
The cut ideal $I_G$ is generated in degree $\leq 4$ if and only if $G$ is free of $K_5$ minors.
\end{conjecture}

In algebraic statistics, minimal generators of toric ideals
are called {\em Markov bases} \cite{Develin2003, Diaconis1998, Takemura2005}.
Thus, what  Conjectures \ref{conj:minors},
\ref{conj:delete}, \ref{conj:degree2} and
\ref{conj:K_5} are about is the complexity of 
Markov bases for moves among the $\mathbb{N}$-valued 
functions on the cuts of a graph $G$. As we shall see 
in Sections 4 and 5, the underlying toric models 
\cite[\S 1.2]{Pachter2005} are important in statistics,
and this endows  our computations and conjectures in this section
with an applied relevance.

From the more theoretical perspective of commutative algebra, 
it appears that Conjecture \ref{conj:K_5} also captures the
class of graphs having normal and Cohen-Macaulay cut ideals.

\begin{conjecture}\label{conj:normal}
The semigroup algebra  
$\kk[q]/I_G$ is normal if and only if $\kk[q]/I_G$ is Cohen-Macaulay if and only if $G$ is free of $K_5$ minors.
\end{conjecture}

That $\kk[q]/I_{K_5}$ is not normal and not Cohen-Macaulay
can be seen in Table 1. The gap between the codimension $(5)$
and the projective dimension $(15)$ is remarkably large in 
this case (we note that the associated semigroup $Q$ and its saturation $Q_{\rm sat}$ differ by only one point). The property of being normal 
is preserved when passing from a semigroup algebra to
a facial subalgebra. Hence we can deduce from
Lemma \ref{lemma:induced} that 
every graph with a $K_5$ minor has a non-normal cut ring $\kk[q]/I_G$.  
Thus, to prove a large part of Conjecture \ref{conj:normal} it would be sufficient to prove 
that graphs $G$  which are free of $K_5$ minors have normal 
semigroup algebras $\kk[q]/I_G$.   Here we are using Hochster's 
Theorem, which states that
normal implies Cohen-Macaulay among semigroup algebras \cite{Hochster1972}.

One question that remains is to characterize 
those $K_5$-free graphs $G$ whose
cut ideal $I_G$ is Gorenstein. Being Gorenstein
seems to depend in a complicated way on the structure of
the graph $G$. In general,
the Gorenstein property is not
preserved under taking toric fiber products and, in particular, 
is not preserved under taking clique sums of graphs.  We do not have a 
firm conjecture on the structure of 
those graphs whose cut ideal is Gorenstein.


\section{From Cut Ideals to Binary Graph Models}\label{sec:graph}

We now explain the correspondence between certain 
cut ideals and the toric ideals of binary graph models.
These are statistical models for $ 2 \times 2 \times 
\cdots \times 2$-contingency tables, whose
algebraic properties were studied by
Develin and Sullivant in \cite{Develin2003}. Our main result
in this section (Theorem \ref{thm:cut2graph}) states
that binary graph models on $n$ nodes
coincide with cut ideals of those graphs 
 on $n+1$ nodes
where one node is connected to all others.

Let $G$ be a graph with vertex set $V = [n] = \{1,2,\ldots,n\}$ 
and edge set $E$, and suppose that $G$ has no isolated vertices.
We introduce a polynomial ring with $2^n$ unknowns,
$$\kk[\,p\,] \,\,\,=\,\,\, 
\kk[\,p_{i_1 i_2 \cdots i_n} \, \, | \, \, i_1,i_2,\ldots,i_n \in \{0,1\}
\,], $$
and a polynomial ring with $\,4 \cdot |E|\,$ unknowns,
$$\kk[\,b\,] \,\,\, =\,\,\,
 \kk[\, b^{e}_{ij}    \, \, | \, \, i,j  \in \{0,1\},  e \in E\,].$$
The binary graph model  is defined by the 
following homomorphism of polynomial rings:
$$\psi_G \,:\, \kk[p] \rightarrow  \kk[b]\,, \quad
p_{i_1 \cdots i_n} \mapsto  \prod_{\{k,l\} \in E } b^{kl}_{i_k i_l}. $$
The kernel of $\psi_G$ is a toric ideal which we denote by $J_G$.
The {\em binary graph model} of $G$ is the zero set of
$J_G$ in $\,\PP^{2^n-1}$. In statistics, this
toric variety corresponds to the
hierarchical model for  $2 \times 2 \times \cdots \times 2$ 
contingency tables where the $2 \times 2$-margins
on the edges of $G$ are fixed.
The Markov basis for this model consists of
the minimal generators of $J_G$.

The {\em suspension} of the graph 
$ G = (V,E)$ is the new graph
 $\widehat{G} $ whose vertex set equals 
$\,[n+1] = V \cup \{n+1\}\,$ and 
whose edge set equals $\,E \cup \{ \{i,n+1\} \, | \, i \in V \}$.
Given any binary string 
$\,{\bf i} = i_1 i_2 \cdots i_n \in \{0,1\}^n$, we
define the associated partition $\,A({\bf i}) | B({\bf i})\,$ 
of  $[n+1]$ by
the condition $k \in B({\bf i})$ if and only if $i_k = 1$.  Similarly, 
if $A|B$ is a partition of $[n+1]$, with $n+1 \in A$,
we get a binary string ${\bf i}(A|B)$ by reversing this  procedure.
This specifies a natural bijection
between the $2^n$ unknowns $p_{\bf i}$ in $\,\kk[p] \,$ and the
$2^n$ unknowns $q_{A|B}$ in $\,\kk[q] $.

\begin{theorem} \label{thm:cut2graph}
Let $\gamma$ be the  ring isomorphism $\, \kk[p] \rightarrow \kk[q]$ 
defined by $\,p_{\bf i}  \mapsto  q_{A({\bf i}) | B({\bf i})}$.  Then,  
$$\gamma(J_G) \quad = \quad  I_{\widehat{G}}.$$
\end{theorem}

We note that this theorem is already known at the level of the underlying
convex polytopes. This is the content of Chapter 5 of \cite{Deza1997}.
The polytope underlying the toric ideal $J_G$ is the
{\em marginal polytope} or {\em covariance polytope} of the graph $G$.
It is isomorphic to the cut polytope of the suspension $\widehat{G}$
under the {\em covariance mapping}, as explained in
\cite[\S 5.2]{Deza1997}. The identification of
$J_G$ with $I_{\widehat{G}}$ in Theorem \ref{thm:cut2graph}
lifts the covariance mapping to the setting of
toric algebra.  Before presenting the proof, we 
discuss a few examples.

\begin{example}
Let $G = K_3$ be the complete graph on three nodes.
The homomorphism $\psi_G$ takes the polynomial ring
$\,\kk[\, p_{000}, p_{001}, p_{010}, p_{011}, 
p_{100}, p_{101}, p_{110}, p_{111}\,]\,$ to the
polynomial ring $\, \kk[\,
b^{12}_{00}, b^{12}_{01}, b^{12}_{10}, b^{12}_{11},
b^{13}_{00}, b^{13}_{01}, b^{13}_{10}, b^{13}_{11},
b^{23}_{00}, b^{23}_{01}, b^{23}_{10}, b^{23}_{11}\,]\,$
by sending
$\, p_{ijk}\,$ to $\,  b^{12}_{ij} b^{13}_{ik} b^{23}_{kl} $.
The kernel $J_G$ is the principal ideal generated by
$\,p_{000} p_{011} p_{101} p_{110} - p_{001} p_{010} p_{100} p_{111}$.
The isomorphism $\gamma$ sends
$\, p_{000} \mapsto q_{1234|},\,
p_{001} \mapsto q_{124|3},\,
p_{010} \mapsto q_{134|2},\,
p_{011} \mapsto q_{14|23},\,
p_{100} \mapsto q_{1|234},\,
p_{101} \mapsto q_{13|24},\,
p_{110} \mapsto q_{12|34} , \,
p_{111} \mapsto q_{123|4}$.
The image of $J_{K_3} $ under $\gamma$ is the
principal ideal $I_{K_4}$ which is discussed in
Example \ref{K4}. Note that $K_4 $ is the suspension of $K_3$. \qed
\end{example}

\begin{example}
Theorem \ref{thm:cut2graph} explains some of the coincidences between rows in our 
Table 1  and the table on page 447 of \cite[\S 2]{Develin2003}. For instance,
the ideal $J_{K_4} \cong I_{K_5} $ is minimally generated by  $20$ quartics
and $40$ sextics. Or, if $G$ is the edge graph of the bipyramid, denoted
$BP$ in \cite{Develin2003}, then its suspension
$\widehat{G} $ is the graph $ G_{207}$ in our Table 1, and
the ideal $\,J_{BP} \cong I_{G_{207}}\,$ is minimally generated by
eight quadrics, $436$ quartics and $2872$ sextics. \qed
\end{example}

The results in Section 3 of \cite{Develin2003} imply the following corollary for cut ideals.
Note that it is consistent with Conjecture \ref{conj:K_5}
because the relevant suspensions $\widehat{G}$ have no $K_5$ minors.

\begin{corollary} \label{fromDevelin}
Let $G$ be a cycle $C_n$ or a complete bipartite graph $K_{2,n}$.  Then the cut ideal $I_{\widehat{G}}$ of the suspension $\widehat{G}$ is generated by binomials having degrees $2$ and $4$.
\end{corollary}

The results in Section 4 of \cite{Develin2003} provide counterexamples to a conjecture that seems to be implied by Table 1; namely, there exist graphs whose cut ideals have minimal generators of odd degree.  The smallest such example for a binary graph model concerns the graph $G = K_2 \times K_3$, the edge graph of the triangular prism, whose graph ideal $J_G$ has a minimal generator of degree $3$.  The suspension of this graph, which has seven vertices, has a cut ideal with an odd degree minimal generator.

\smallskip 

\begin{proof}[Proof of Theorem \ref{thm:cut2graph}]
It suffices to show that there are a pair of homomorphisms $\alpha: \kk[b] \rightarrow \kk[s,t]$ and $\beta: \kk[s,t] \rightarrow \kk[b]$ such that $\phi_{\widehat{G}} \circ \gamma = \alpha \circ \psi_{G}$ and 
$\psi_G \circ \gamma^{-1} =  \beta \circ \phi_{\widehat{G}}$. 
The maps $\alpha$ and $\beta$, restricted to $\kk[p]/J_G$ and $\kk[q]/I_{\widehat{G}}$ respectively, will then lift to  the isomorphism $\gamma$.  To do this correctly,
 we extend $\kk[s,t]$ and $\kk[b]$ to allow fractional powers of the unknowns.
 Which fractional powers are needed will be clear from the context.

We define the map
$\,\alpha:  \kk[b] \rightarrow \kk[s,t]\,$ as follows:
$$b^{kl}_{00}  \mapsto t_{kl} t^{\frac{1}{\deg(k)}}_{k,n+1} t^{\frac{1}{\deg(l)}}_{l, n+1}, 
 \quad  b^{kl}_{01} \mapsto    s_{kl} t^{\frac{1}{\deg(k)}}_{k,n+1} s^{\frac{1}{\deg(l)}}_{l, n+1}$$
 $$b^{kl}_{10}  \mapsto s_{kl} s^{\frac{1}{\deg(k)}}_{k,n+1} t^{\frac{1}{\deg(l)}}_{l, n+1}, 
 \quad  b^{kl}_{11} \mapsto    t_{kl} s^{\frac{1}{\deg(k)}}_{k,n+1} s^{\frac{1}{\deg(l)}}_{l, n+1}.$$
 Here ${\rm deg}(k)$ denotes the degree of the node $k$ in the
 graph $G$, and similarly for the node~$l$.
 
We wish to show that $\alpha$ satisfies 
$\phi_{\widehat{G}} \circ \gamma = \alpha \circ \psi_{G}$.
To do this, we look at which unknowns $s_{kl}, t_{kl}$ appear to which powers in the monomials
$\,\alpha (\psi_{G} ( p_{\bf i}))\,$ and $\, \phi_{\widehat{G}} ( \gamma (p_{\bf i}))$.
An unknown $s_{kl}$ appears in $\,\alpha ( \psi_G ( p_{\bf i} ))\,$ 
with multiplicity one if and only if $\,i_k i_l \in \{ 01,10\}\,$
 if and only if $\{k,l\} \in {\rm Cut}(A({\bf i}) | B({\bf i}) )\,$
  if and only if $s_{kl}$ appears in $\, \phi_{\widehat{G}}( \gamma (p_{\bf i}))\,$ with multiplicity one.  
  A similar argument shows that $t_{kl}$ appears with the same multiplicity in both 
  $\,\alpha (\psi_G ( p_{\bf i} ))\,$ and $\, \phi_{\widehat{G}} ( \gamma (p_{\bf i}))$.  
  To check the multiplicity of $s_{k,n+1}$ (and similarly for $t_{k,n+1}$), note that the fractional powers guarantee that $s_{k,n+1}$ appears in   $\,\alpha( \psi_G ( p_{\bf i} ))\,$ 
  if and only if it has multiplicity one in $\,\alpha( \psi_G ( p_{\bf i} ))$. This happens
  if and only if $i_k = 1$ if and only if $(k,n+1) \in {\rm Cut}(A({\bf i}) | B({\bf i}) )$ if and only if $s_{kl}$ appears in $ \,\phi_{\widehat{G}} ( \gamma (p_{\bf i}))\,$ with multiplicity one.

 We now define our second ring homomorphism   $\,\beta: \kk[s,t] \rightarrow \kk[b]\,$ as follows:
 $$
 s_{k,n+1} \mapsto  \prod_{l : \{k,l\} \in E }  \left( b^{kl}_{00}  b^{kl}_{01} \right)^{-\frac{1}{2} } \cdot B, \quad \quad 
 t_{k,n+1} \mapsto  \prod_{l : \{k,l\} \in E }  \left( b^{kl}_{10}  b^{kl}_{11} \right)^{-\frac{1}{2} } \cdot B
 $$ 
 $$
s_{kl}  \mapsto \left( b^{kl}_{01}  b^{kl}_{10} \right)^{\frac{1}{2}}, \quad \quad  
t_{kl}  \mapsto \left( b^{kl}_{00}  b^{kl}_{11} \right)^{\frac{1}{2}}. $$
Here $B$ denotes the product of all unknowns in $\kk[b]$ raised to the power $1/2n$:
$$B =  \prod_{\{k,l\} \in E}  \prod_{i,j \in \{0,1\}}  (b^{kl}_{ij} )^{\frac{1}{2n}}. $$
To prove that $\beta$ satisfies 
$\psi_G \circ \gamma^{-1} =  \beta \circ \phi_{\widehat{G}}$ we compare the multiplicity of $b^{kl}_{ij}$ in $\,\psi_G ( \gamma^{-1}(q_{A|B}))\,$ and  $\,\beta ( \phi_{\widehat{G}}(q_{A|B}))$. 
By symmetry, it suffices to analyze the case $ij = 00$. For fixed $k,l$, the 
unknown $b^{kl}_{00}$ has multiplicity one in 
$\,\psi_G (\gamma^{-1}(q_{A|B}))\,$ if and only if $\{k,l\} \notin {\rm Cut}(A|B)$ and $k,l \in A$.
Here, $b^{kl}_{01}$, $b^{kl}_{10}$, $b^{kl}_{11}$ all occur with multiplicity zero.  

Now we analyze the multiplicity of $b^{kl}_{ij}$ in $\,\beta(\phi_{\widehat{G}}(q_{A|B}))$.  Suppose we are in the case $\{k,l\} \notin {\rm Cut}(A|B)$ and $k,l \in A$.  This means that $t_{kl}t_{k,n+1}t_{l,n+1}$ is a factor of $\phi_{\widehat{G}}(q_{A|B})$.  Looking at the expansion of $\beta( \phi_{\widehat{G}} (q_{A|B}))$, aside from the factor $B^n$, the only multiplicands which possibly contain $b^{kl}_{00}$ are $t_{kl}, t_{k,n+1}, t_{l,n+1}$.  The first contributes 
$(b^{kl}_{00})^{\frac{1}{2}}$, the second and third contribute nothing, and the factor of $B^n$ contribute  $(b^{kl}_{00})^{\frac{1}{2}}$ for a grand total of $b^{kl}_{00}$. 
 On the other hand, $b^{kl}_{01}$ appears with multiplicity zero because $t_{kl}$ and $t_{k,n+1}$  contribute nothing, $t_{l,n+1}$ contributes $(b^{kl}_{01})^{-\frac{1}{2}}$,
 and $B^n$ contributes $(b^{kl}_{00})^{\frac{1}{2}}$.  A similar argument shows that $b^{kl}_{10}$ and $b^{kl}_{11}$ also appear with multiplicity zero.  This agrees with the multiplicity of $b^{kl}_{ij}$ in $\,\psi_G ( \gamma^{-1}(q_{A|B}))$. This completes the proof of Theorem \ref{thm:cut2graph}.
\end{proof}


\section{From Jukes-Cantor Phylogenetic Models to Cut Ideals}\label{sec:split}

In this section we apply cut ideals to phylogenetics.
Our main result (Theorem \ref{thm:cutsplit}) states that cut ideals 
of graphs with $n$ nodes are precisely the
binary Jukes-Cantor models on cyclic
split systems on $n$ taxa. This class includes the Jukes-Cantor models on
phylogenetic trees whose algebraic properties were studied in
\cite{BW} and  \cite{Sturmfels2005}. We rederive
the quadratic Gr\"obner basis for these ideals
by relating Theorem \ref{thm:cutfiber} to
 \cite[Theorem 21]{Sturmfels2005}.

The extension of statistical models of evolution from phylogenetic trees
to split systems is due to David Bryant, who described these models in
\cite{Bryant2005}. This extension has the double advantage of
being useful for biological applications and leading to
a richer mathematical theory. In what follows we give
an algebraic introduction to Jukes-Cantor models for 
arbitrary split systems. Later on we specialize
to split systems which are cyclic, and hence most
relevant for the {\tt NeighborNet}  method
\cite{Bryant2004}. This will take us back to
cut ideals.

\subsection{The one-parameter model associated with a single split}

We consider a set of $n$ taxa labeled by $[n] = \{1,2,\ldots,n\}$.
Each Jukes-Cantor model is a subvariety of the
$(2^{n}-1)$-dimensional projective space $\PP^{2^{n}-1}$ whose
coordinates we denote by $\,p_{i_1 \cdots i_n}$.
The coordinate $\,p_{i_1 \cdots i_n}\,$ represents the probability 
of observing the states $i_1,\ldots,i_n \in \{0,1\}$
at the taxa. We shall employ a linear change of
coordinates which is known as the {\em Fourier transform}
or {\em Hadamard conjugation}; see \cite[\S 4.4]{Pachter2005} and \cite[\S 2]{Sturmfels2005}.
The Fourier coordinates are here denoted
$f_{j_1 \ldots j_n}$, and they are related to the 
probability coordinates as follows:
\begin{equation}
\label{qTOp}
f_{j_1 \cdots j_n} \quad = \quad
\sum (-1)^{i_1 j_1 + \cdots + i_n j_n } \cdot p_{i_1 \cdots i_n},
\end{equation}
where the sum is over all elements $(i_1,\ldots,i_n)$
of the abelian group $\, (\zz/2\zz)^{n}$.
It is very easy to invert this linear transformation. Namely, we have
\begin{equation}
\label{pTOq}
 p_{i_1 \cdots i_n} \quad  = \quad \frac{1}{2^{n}}
 \sum (-1)^{j_1 i_1 + \cdots + j_n i_n } \cdot f_{j_1 \cdots j_n},
\end{equation}
where the sum is over $\,(j_1,\ldots,j_n) \in (\zz/2\zz)^{n}$.

A {\em split} $\{C,D\}$ is a partition 
$\, C \cup D \, = \, \{1,\ldots,n\}$ of the set 
of taxa such that $\, n \in D$. We fix
a split $\{C,D\}$ and we introduce one free parameter $u$. In statistical
applications, this parameter $u$ would range over real numbers
between $0$ and $\frac{1}{2}$. In algebraic geometry we
allow any point $(u_0:u_1)$ on  the
complex projective line $\PP^1$,
with $u_0 = 1$ and $u_1 = u$.

We map the $u$-line  $\PP^1$ into
the probability space $\PP^{2^{n}-1}$ by setting
\begin{itemize}
\item
$\, f_{j_1 \cdots j_n} \quad = \quad 0\quad \,\,$ if $\,  j_1 + \cdots + j_n \,$ is odd.
\item
$\, f_{j_1 \cdots j_n} \quad = \quad u_0\quad$ if $\,\sum_{k \in C} j_k \,$  and
$\,\sum_{k \in D} j_k \,$  are both even.
\item
$\, f_{j_1 \cdots j_n} \quad = \quad u_1 \quad $ if $\,\sum_{k \in C} j_k \,$  and
$\,\sum_{k \in D} j_k \,$  are both odd.
\end{itemize}
This line in $\PP^{2^{n}-1}$ is the Jukes-Cantor model
associated with the split $\{C,D\}$. Using the  transformation (\ref{qTOp}),
we can express the parameterization
in probability coordinates:
\begin{itemize}
\item $\,p_{i_1 \cdots i_n} \,\, = \,\, (u_0+u_1)/4 \quad $ if
$\, i_1 = \cdots = i_n $.
\item  $\,p_{i_1 \cdots i_n} \,\, = \,\, (u_0-u_1) /4\quad $ if
$\,i_k = 1\,$ for all $k \in C\,$ and 
$\,i_l = 0\,$ for all $l \in D $.
\item  $\,p_{i_1 \cdots i_n} \,\, = \,\, (u_0-u_1) /4\quad $ if
$\,i_k = 0\,$ for all $k \in C\,$ and 
$\,i_l = 1\,$ for all $l \in D $.
\item  $\,p_{i_1 \cdots i_n} \,\, = \,\, 0 \quad $ in all other cases.
\end{itemize}
In summary, the Jukes-Cantor model for a single split is a straight line
in $\PP^{2^n-1}$. Given two points in this model, 
we can multiply their Fourier coordinates,
one coordinate at a time, and we get a new point in the model.
Thus the model is a semigroup with respect to multiplication
of Fourier coordinates. The model is a line which is also a toric curve.

\subsection{The Jukes-Cantor model defined by an arbitrary split system}

A {\em split system} is simply a collection of 
$r$ distinct splits of $[n] = \{1,\ldots,n\}$, for some positive integer $r$:
$$ \Sigma \quad = \quad \bigl\{
\,\{C_1,D_1\}, \{C_2,D_2\},\ldots, \{C_r,D_r\} \bigr\}. $$
Each split $\{C_i,D_i\}$ specifies a one-parameter Jukes-Cantor
model, which is  a semigroup under multiplication of Fourier coordinates.
We define the {\em Jukes-Cantor model of $\Sigma$}  to be the semigroup
generated by the $r$ one-parameter models of the splits $\{C_i,D_i\}
\in \Sigma $.

Explicitly, the parametrization of this Jukes-Cantor model 
is given as follows. The parameter space is the direct product of
$r$ copies of the projective line $\PP^1$. The homogeneous coordinates
of the $i$-th projective line $\PP^1$ are denoted $(u^i_{0}:u^i_1)$.
There are precisely  $2^{n-1}$ nonzero Fourier coordinates 
$\,f_{j_1 \cdots j_n}$. They are indexed by the group
$$ (\zz/2\zz)^{n}_{\rm even} \quad = \quad
\bigl\{\, (j_1,\ldots,j_n)\in (\zz/2\zz)^{n} \,:\,
j_1 + \cdots + j_n \,\,\hbox{is even} \,\bigr\}. $$
Each nonzero Fourier coordinate is expressed as
a monomial of degree $r$ in the parameters:
\begin{equation}
\label{jkPara}
f_{j_1 \ldots j_n} \quad = \quad
 \prod_{\{C_i,D_i\} \in \Sigma}  u^i_{\sum_{k \in C_i } j_k} .
\end{equation}
Since this parametrization is given by monomials, the ideal of
algebraic invariants of the Jukes-Cantor model is a toric ideal
in the Fourier coordinates.
This toric ideal is the kernel of the ring map (\ref{jkPara})
and we denote it by $JC_\Sigma$.
It lives in the polynomial ring $\kk[\,f \,]$ whose
generators are the $2^{n-1}$ Fourier coordinates
$\,f_{j_1 \cdots j_n}$  indexed by
$\, (\zz/2\zz)^{n}_{\rm even}$.

It is important to understand that Jukes-Cantor models
are toric varieties, since $JC_\Sigma$ is a toric ideal
in the Fourier coordinates, but Jukes-Cantor models are \emph{not}
toric models (i.e. log-linear models or discrete exponential families) 
in the sense of \cite[\S 1.2]{Pachter2005} because
$JC_\Sigma$ is not a toric ideal when rewritten in the probability coordinates 
$p_{i_1 \ldots i_n}$ via the Fourier transform (\ref{qTOp}).

\begin{proposition}
If $\Sigma$ consists of $r$ splits then
the Jukes-Cantor model is $r$-dimensional.
\end{proposition}

\begin{proof}
We can write the $2^{n-1}$ nonzero monomials in the parametrization
(\ref{jkPara}) as the columns of a zero-one matrix $A$
with $2r$ rows, one for each unknown $u^i_0$ 
and $u^i_1$, as in \cite{Sturmfels1995}
or in \cite[\S 1.2]{Pachter2005}. The rows of this matrix
span an $r+1$-dimensional linear space.   This implies that
the semigroup algebra $\,\kk[f]/JC_\Sigma\,$ has Krull dimension $r+1$,
and hence the associated projective variety
(which is our Jukes-Cantor model) has dimension $r$.
\end{proof}

Jukes-Cantor models for split systems
do indeed generalize the familiar models associated with trees.
Let $T$ be a tree with leaves labeled by $[n]$.
Every edge of $T$ defines a split $\{C,D\}$ of $[n]$. We write
$\Sigma(T)$ for the set of splits coming from all the edges of $T$.

\begin{remark} 
$JC_{\Sigma(T)}$ equals the usual Jukes-Cantor model
associated with the tree $T$.
\end{remark}

\begin{proof}
This is seen by
comparing the parametrization for split systems in (\ref{jkPara}) with that
given in \cite[\S 3]{Sturmfels2005} for group based models on trees.
The condition that  $\sum_{k \in C_i } j_k$ is even in the split system
representation is replaced with the condition that $\sum_{k \in
  \Lambda(e)} j_k$ is even where $\Lambda(e)$ is the set of leaves
below the edge $e$.  The concept of being a ``leaf below an edge'' is
equivalent to being on one side of a split.
\end{proof}


\subsection{Cyclic Split Systems}

We now turn our attention to the family of cyclic split systems. These
split systems are particularly useful for representing and analyzing 
metric spaces in biology,  as
they can be drawn in the plane using {\tt NeighborNet} \cite{Bryant2004}.

Formally, we define cyclic split systems as follows.
We draw a convex $n$-gon in the plane and label the
vertices by $1,\ldots,n$ in clockwise order.
Every line in the plane that does not pass through any of the
vertices defines a split $\{C,D\}$. The {\em complete cyclic
split system} $\Sigma^{(n)}$ is the collection of all
splits of $[n] = \{1,\ldots,n\}$ which arise in this manner. 

\begin{remark} The number of non-trivial cyclic splits
in $\Sigma^{(n)}$ equals 
$\,n(n-1)/2 $.
\end{remark}

A {\em cyclic split system} is any subset of $\Sigma^{(n)}$.
 In other words, a split system
$\Sigma$ is cyclic if, for each split
$\{C,D\} \in \Sigma$, the set $C$ is an interval of integers
$C = [k,l] = \{k, k+1, \ldots, l \}$.

Now we will show that every cyclic split ideal 
$JC_\Sigma$ is a cut ideal.
We associate with each cyclic split system
 $\Sigma$ a graph $G_\Sigma$ with vertex
set $[n]$ as follows. For each cyclic split $\{C,D\} \in \Sigma$
where $n \in D$ and $C = [k,l]$, we introduce the edge $\{k-1,l\}$ in $G_\Sigma$.
(Here $0 := n$). Thus $G_\Sigma$ is a graph with one edge for each split in $\Sigma$.
The representation of a cyclic split system $\Sigma$ by its graph
$G_\Sigma$ is very natural as the following proposition shows.

\begin{proposition} \label{treesubd}
Let $T$ be a planar tree with leaves 
labeled cyclically $1,\ldots,n$ 
and  $\Sigma(T)$
the associated cyclic split system.  Then 
the graph $G_{\Sigma(T)}$ consists of the edges
in the subdivision of the convex $n$-gon 
which is dual to  the tree $T$. 
\end{proposition} 

\begin{proof}
The proof of this result is straightforward. The idea is 
illustrated in Figure 1.
\end{proof}

\begin{figure}[h] \label{treefig}
\begin{center}
\includegraphics{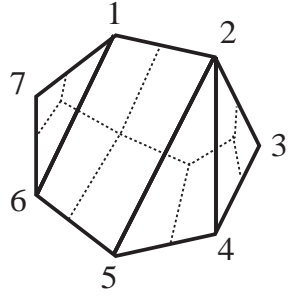}
\caption{A tree with seven leaves and corresponding subdivision of the 7-gon}
\end{center}
\end{figure}

We now come to the main result in this section.
We define a bijection between the set of all
$2^{n-1}$ cuts  of $[n]$ and the
set $(\zz/2\zz)^{n}_{\rm even}$ of binary strings that sum to zero.
If $A|B$ is any cut then the corresponding binary string
$j_1 j_2 \cdots j_n$ is defined as follows:
$$
\hbox{$ j_k = 1\,\,$  if \,\, $\{k-1,k\} \in {\rm Cut}(A|B)\,\,$
and $ \,\,j_k = 0 \,$  otherwise} $$
It is easy to see that $\,j_1 + \cdots + j_n$ is even and that the cut $A|B$ is uniquely encoded in
the string $j_1 j_2 \cdots j_n$. This bijection
defines an isomorphism of polynomial rings
$\,\tau : \kk [ q ] \rightarrow \kk[f] \,$
by sending the unknown $\,q_{A|B}\,$ to $\,f_{j_1 \cdots j_n}$.

\begin{theorem} \label{thm:cutsplit}
Let $\Sigma$ be a cyclic split system and $G_\Sigma $ the associated
graph.  Then the Jukes-Cantor model  $JC_\Sigma$
equals the image of the cut ideal 
$I_{G_\Sigma}$ under the isomorphism
$\,\tau$.
\end{theorem}

\begin{proof} 
To see that the preceding bijection between cut coordinates and
Fourier coordinates gives an isomorphism between the cut model for
$G_\Sigma$ and the Jukes-Cantor model $JC_\Sigma$ we must define an appropriate bijection between the parameters.  This
bijection between parameters is induced by the map which sends a
cyclic split in $\Sigma$ to an edge of the graph $G_\Sigma$.  Namely,
we identify the $\mathbb{P}^1$  parameter space associated to the
split $\{C_i,D_i\}$, 
where $C_i = 
[k+1,l]$, with the $\mathbb{P}^1$ parameter space associated to the
edge $\{k,l\}$ in $G_\Sigma$ via
\begin{equation}
\label{iiddii}
 (u^i_{0}:u^i_{1})\,\, =\,\, (t_{kl}: s_{kl}). 
\end{equation}
Now, the unknown $s_{kl}$ appears in  the squarefree monomial
$\,\phi_{G_\Sigma}(q_{A|B})\,$ if and only if
$\,\{k,l\} \in {\rm Cut}(A|B)\,$ if and only if 
$\,j_{k+1} + \cdots + j_l \,$ is odd if and only if 
$\,\sum_{\nu \in C_i} j_\nu \,$ is odd if and only if
the unknown $u^i_1$ appears in  the squarefree monomial
on the right hand side of (\ref{jkPara}).
Likewise, $t_{kl}$ appears in
$\,\phi_{G_\Sigma}(q_{A|B})\,$ if and only if
$u^i_0$ appears in  the right hand side of (\ref{jkPara}).
This shows, modulo the identification (\ref{iiddii}),
that the image of the cut coordinate
$q_{A|B}$ under the map $\phi_{G_\Sigma}$
equals the image of the Fourier coordinate
 $f_{j_1 \cdots j_n}$ under the map (\ref{jkPara}).
Therefore, both maps have the same kernel, and
we conclude $\,JC_\Sigma  = \tau( I_{G_\Sigma})$.
\end{proof}

\begin{example} \label{ex:k4cut}
Let $\Sigma = \Sigma^{(4)}$ be the complete
cyclic split system on four taxa, i.e.,
\begin{equation}
\label{splitorder}
\Sigma \,\,\, = \,\,\,
\bigl\{\{12,34\},\{23,14\}, \{1,234\}, 
\{2,134\}, \{3,124\} ,\{123,4\} \bigr\}. \end{equation}
The  associated graph $G_\Sigma$ is the complete graph on $\{1,2,3,4\}$.
With the ordering of the splits as in (\ref{splitorder}),  the map $\tau$ and the
 Jukes-Cantor parametrization (\ref{jkPara}) are given by
\begin{eqnarray*}
q_{|1234} & \mapsto \quad
f_{0000} \,\,\, & \mapsto \quad  u^4_{0} \cdot  u^2_{0} \cdot  u^3_{0} \cdot  u^5_{0} \cdot  u^1_{0} \cdot  u^6_{0} \\
q_{4|123} & \mapsto \quad
f_{1001} \,\,\, & \mapsto \quad  u^4_{0} \cdot  u^2_{0} \cdot  u^3_{1} \cdot  u^5_{0} \cdot  u^1_{1} \cdot  u^6_{1} \\
q_{3|124} & \mapsto \quad
f_{0011} \,\,\, & \mapsto \quad  u^4_{0} \cdot  u^2_{1} \cdot  u^3_{0} \cdot  u^5_{1} \cdot  u^1_{0} \cdot  u^6_{1} \\
q_{2|134} & \mapsto \quad
f_{0110} \,\,\, & \mapsto \quad  u^4_{1} \cdot  u^2_{0} \cdot  u^3_{0} \cdot  u^5_{1} \cdot  u^1_{1} \cdot  u^6_{0} \\
q_{1|234} & \mapsto \quad
f_{1100} \,\,\, & \mapsto \quad  u^4_{1} \cdot  u^2_{1} \cdot  u^3_{1} \cdot  u^5_{0} \cdot  u^1_{0} \cdot  u^6_{0} \\
q_{12|34} & \mapsto \quad
f_{1010} \,\,\, & \mapsto \quad  u^4_{0} \cdot  u^2_{1} \cdot  u^3_{1} \cdot  u^5_{1} \cdot  u^1_{1} \cdot  u^6_{0} \\
q_{13|24} & \mapsto \quad
f_{1111} \,\,\, & \mapsto \quad  u^4_{1} \cdot  u^2_{0} \cdot  u^3_{1} \cdot  u^5_{1} \cdot  u^1_{0} \cdot  u^6_{1} \\
q_{14|23} & \mapsto \quad
f_{0101} \,\,\, & \mapsto \quad  u^4_{1} \cdot  u^2_{1} \cdot  u^3_{0} \cdot  u^5_{0} \cdot  u^1_{1} \cdot  u^6_{1}. 
\end{eqnarray*}
Under the identification  (\ref{iiddii}), this coincides with the
parametrization in Example \ref{K4}. 
The Jukes-Cantor ideal 
for the complete split system on four taxa equals
 $$ J_\Sigma \quad = \quad \langle \,
  f_{0000} f_{0101} f_{1010} f_{1111} 
 \,-\,   f_{0011} f_{0110} f_{1001} f_{1100} \,\rangle . $$
 The ordering of the factors $u^i_j$ in the above monomials
coincides with the lexicographic ordering 
of the edges of $K_4$. If we set 
$\,u^1_0 = u^1_1 = u^2_0 = u^2_1 = 1 \,$
in the parametrization, then we
get the $4$-cycle in Example \ref{C4}, which represents the
Jukes-Cantor model for the star tree.
This model is the same as the rooted
claw tree $K_{1,3}$ in 
\cite[Example 14]{Hosten2005}. \qed
\end{example}

\subsection{Algebraic invariants for Jukes-Cantor models on
cyclic split systems}

The polynomials in the ideal $J_\Sigma$ are known as
{\em algebraic invariants} in phylogenetics. When
expressed in terms of the coordinates $p_{i_1 \cdots i_n}$
via (\ref{qTOp}), these polynomials are the algebraic
relationships which hold among the joint probabilities
for all distributions in the model. 
Using Theorem  \ref{thm:cutsplit}, we can now translate
our results and conjectures about cut ideals to the
setting of Jukes-Cantor models.
We begin by giving a new proof of a known result.

\begin{corollary} {\rm \cite[Theorem 2 (a)]{Sturmfels2005}}
Consider the Jukes-Cantor model
for any trivalent tree $T$ with taxa $[n]$. Then
the ideal $JC_{\Sigma(T)}$ has
a Gr\"obner basis consisting of quadrics.
\end{corollary}

\begin{proof}
By Proposition \ref{treesubd} and Theorem \ref{thm:cutsplit},
we have $JC_{\Sigma(T)} = I_G$,
where $G$ is the edge graph of 
a triangulation of the $n$-gon.
Such a planar graph can be 
decomposed into triangles
using $2$-sums. The result hence
follows from Theorem \ref{thm:cutfiber}.
\end{proof}

We now discuss the Jukes-Cantor ideals $JC_\Sigma$
for some other cyclic split systems. Each of the graphs $G$
in Table 1 corresponds to such a split system.
Namely, for each edge $\{k,l\}$ of $G$
we introduce the cyclic split $\{C,D\}$
where $C = \{k+1,k+2,\ldots,l\}$ and
$D = [n] \backslash C$.

The complete graph $K_n$ corresponds to the complete
split system $\Sigma^{(n)}$. Table 1 reveals that
the algebraic invariants for $\Sigma^{(5)}$ are generated
in degree $\leq 6$ and the algebraic invariants for $\Sigma^{(6)}$
are generated in degree $\leq 10$.
Conjectures \ref{conj:degree2} and Conjecture
\ref{conj:K_5} translate into conjectures
as to which Jukes-Cantor ideals $JC_\Sigma$
are generated by quadrics and which are generated by
quartics. Whether the generating degree $\mu(JC_\Sigma)$
for a cyclic split system
can only decrease upon removal of a split is still unknown,
in light of Conjecture \ref{conj:delete}.

Huson and Bryant \cite{Huson2006} have shown that
cyclic split systems, even if they do not arise from trees, always
have useful representations by  phylogenetic networks. However,
this representation is generally not unique
\cite[Figure 5]{Huson2006}.
These {\em split networks} on $n$ taxa are thus in
many-to-one correspondence, via Theorem \ref{thm:cutsplit}
to graphs with $n$ vertices, and our results here shed light
on the algebraic invariants of the associated statistical 
model \cite{Bryant2005}. One concrete application 
of this correspondence to phylogenetics 
will be the exact computation of
maximum likelihood parameters 
for splits models as described in
\cite[\S 6]{Hosten2005}.

\begin{example}
Let $n=6$ and consider the bipartite graph 
$K_{3,3}$ where the bipartition separates
$\{1,3,5\} $ from $\{2,4,6\}$. The corresponding
split system $\Sigma$ consists of the six
trivial splits $\{\{i\} , [6] \backslash \{i\} \}$
and the three non-trivial splits
$\{123,456\} $, $\{234,156\}$ and $\{ 345,126\}$.
This is the smallest split system 
whose split network is not unique.
It is depicted in \cite[Figure 5]{Huson2006}.
Using our Table 1 in Section 3, we see that the corresponding
Jukes-Cantor ideal $JC_\Sigma$ is minimally
generated by $63$ quadrics and $72$ quartics.
The semigroup algebra $\kk[f]/JC_\Sigma$
is also normal and hence Cohen-Macaulay,
by Hochster's Theorem \cite{Hochster1972}. \qed
\end{example}

\bigskip
\bigskip

\noindent {\bf Acknowledgement:}
Bernd Sturmfels was partially supported by the National Science
Foundation (DMS-0456960). 

\bigskip


\bigskip

\end{document}